\documentclass[12pt]{article}
\usepackage{amsmath, amssymb, amsfonts, amsthm}
\def\CC{{\rm \kern.24em \vrule width.02em height1.4ex
depth-.05ex \kern-.26em C}}

\def\TagOnRight

\def\AA{{\it I}\hskip-3pt{\tt A}}

\def\QQ{\rlap {\raise 0.4ex \hbox{$\scriptscriptstyle |$}}
  {\hskip -0.1em Q}}

\newcommand{\be}{\begin{equation}}
\newcommand{\ee}{\end{equation}}
\newcommand{\bea}{\begin{eqnarray}}
\newcommand{\eea}{\end{eqnarray}}
\newcommand{\Bea}{\begin{eqnarray*}}
\newcommand{\Eea}{\end{eqnarray*}}

\newcommand{\bi}{\begin{itemize}}
\newcommand{\ei}{\end{itemize}}

\newtheorem{Definition}{Definition}[section]
\newtheorem{Theorem}[Definition]{Theorem}
\newtheorem{Lemma}[Definition]{Lemma}

\newtheorem{Corollary}[Definition]{Corollary}

\newtheorem{Revised Question}[Definition]{Revised Question}
\newtheorem{Notation}[Definition]{Notation}

\theoremstyle{remark}

\begin{document}

\title{The mu vector, Morse inequalities and a generalized lower bound theorem for locally tame combinatorial manifolds}
\author{Bhaskar Bagchi\\
Theoretical Statistics and Mathematics Unit\\
Indian Statistical Institute\\
Bangalore - 560 059, India.\\
email: bbagchi@isibang.ac.in}
\date{}

\maketitle

\abstract{In a recent work [2] with Datta, we introduced the mu vector (with respect to a given field) of simplicial complexes and used it to study tightness and lower bounds. In this paper, we modify the definition of mu vectors. With the new definition, most results of [2] become correct without the hypothesis of 2-neighbourliness. In particular, the combinatorial Morse inequalities of [2] are now true of all simplicial complexes.

As an application, we prove the following generalized lower bound theorem (GLBT) for connected locally tame combinatorial manifolds. If $M$ is such a manifold of dimension $d$, then for $1 \leq \ell \leq \frac{d-1}{2}$ and any field $\mathbb{F}, ~ g_{\ell+1} (M) \geq \binom{d+2}{\ell+1} \sum\limits_{i=1}^\ell (-1)^{\ell-i} \beta_i (M;\mathbb{F})$. Equality holds here if and only if $M$ is $\ell$-stacked.

We conjecture that, more generally, this theorem is true of all triangulated connected and closed homology manifolds. A conjecture on the sigma vectors of triangulated homology spheres is proposed, whose validity will imply this GLB Conjecture for homology manifolds. We also prove the GLBC for all connected and closed combinatorial 3-manifolds. Thus, any connected closed combinatorial manifold $M$ of dimension three satisfies $g_2 (M) \geq 10 \beta_1 (M;\mathbb{F})$, with equality iff $M$ is 1-stacked. This result settles a question of Novik and Swartz [6] in the affirmative.}

\section{Introduction and results.}

All simplicial complexes considered in this paper are finite and abstract. All homologies are simplicial homologies with coefficients in an arbitrary field $\mathbb{F}$, which is fixed throughout. In consequence, the Betti numbers (reduced or not), and hence also the sigma- and mu-vectors attached to a simplicial complex, depend on the field $\mathbb{F}$. This dependence will not be indicated in the notations. By a homology manifold (with or without boundary) we shall mean a triangulation of an $\mathbb{F}$-homology compact manifold. The compact homology manifolds without boundary are called the closed homology manifolds.

For any simplicial complex $X,V(X)$ will denote the vertex set of $X$. For any set $A,X[A]$ will denote the subcomplex $\{\alpha \in X : \alpha \subseteq A\}$. Thus $X[A]$ is the induced subcomplex of $X$ with vertex set $V(X) \cap A$. But the extra flexibility, allowed by the possibility that $A$ may not be a subset of $V(X)$, will be useful in the arguments handling improper bistellar moves. For $x \in V(X)$, the link of $x$ in $X$, denoted by $\ell k (x,X)$, is the subcomplex $\{\alpha : x\not\in \alpha, \alpha \sqcup \{x\} \in X\}$.

For any non-empty finite set $\alpha, \overline{\alpha}$ (respectively $\partial \alpha$) will denote the simplicial complex whose faces are all the subsets (respectively all proper subsets) of $\alpha$. Thus, when $\# (\alpha) =d+1, ~ \overline{\alpha}$ is the unique $(d+1)$-vertex triangulation of the $d$-ball (namely, the face complex of a geometric $d$-simplex). When $\# (\alpha) =d+2$, $\partial \alpha$ is the unique $(d+2)$-vertex triangulation of the $d$-sphere (namely, the boundary complex of a geometric $(d+1)$-simplex). The $(d+2)$-vertex $d$-sphere is called the standard $d$-sphere; when we do not wish to indicate its vertex set, it is denoted by $S^d_{d+2}$. We adopt the convention that the empty set is a face (of dimension $-1$) of all simplicial complexes. The join $X_1 \ast X_2$ of two simplicial complexes $X_1$ and $X_2$ is the simplicial complex whose faces are the disjoint unions of the faces of $X_1$ with the faces of $X_2$.

If $\alpha$ and $\beta$ are two disjoint non-empty finite sets with $\# (\alpha) + \#(\beta)=d+2$, then $\overline{\alpha} \ast \partial \beta$ and $\overline{\beta} \ast \partial \alpha$ are two triangulations of the $d$-ball with common boundary $\partial \alpha \ast \partial \beta$. If $X$ is a $d$-dimensional simplicial complex with $\overline{\alpha} \ast \partial \beta$ as an induced subcomplex, then $X$ is said to admit the bistellar move $\alpha \mapsto \beta$. In this case, the simplicial complex $Y=(X\backslash (\overline{\alpha} \ast \partial \beta)) \cup (\overline{\beta} \ast \partial \alpha)$ is said to be obtained from $X$ by the bistellar move $\alpha \mapsto \beta$. Clearly, in this case, $X$ and $Y$ have homomorphic geometric realizations, and $X$ can be obtained from $Y$ by the reverse bistellar move $\beta \mapsto \alpha$. If $\dim (\overline{\beta})=t$, then we say that the bistellar move $\alpha \mapsto \beta$ is of index $t$. Then the reverse move is of index $d-t$. A move of index 0 introduces a new vertex, while a move of index $d$ deletes an old vertex; these are sometimes called the improper bistellar moves.

A combinatorial $d$-sphere is a simplicial complex which may be obtained from the standard $d$-sphere $S^d_{d+2}$ by a finite sequence of bistellar moves. Thus all combinatorial $d$-spheres are triangulations of the topological $d$-sphere. A combinatorial $d$-manifold is a $d$-dimensional simplicial complex all whose vertex links are combinatorial $(d-1)$-spheres. Thus combinatorial $d$-manifolds are triangulations of topological closed $d$-manifolds, though the converse is false.

We now introduce a special class of combinatorial spheres and of combinatorial manifolds.

\begin{Definition} A tame $d$-sphere is a $d$-dimensional simplicial complex which may be obtained from $S^d_{d+2}$ by a finite sequence of bistellar moves, each of index $< d/2$. A locally tame $d$-manifold is a $d$-dimensional simplicial complex all whose vertex links are tame $(d-1)$-spheres.
\end{Definition}

Thus, all tame $d$-spheres are combinatorial $d$-spheres and all locally tame $d$-manifolds are combinatorial $d$-manifolds. Notice that the tame $d$-spheres are just the $\lceil \frac{d}{2}\rceil$-stellated spheres and the tame $d$-manifolds are the members of the class $W_{\lceil \frac{d-1}{2}\rceil} (d)$, as defined in [1, 2]. But, as we shall soon see, there are good reasons to single out these two special classes for detailed study.

Recall that, for $0 \leq k \leq d$ and a $d$-dimensional simplicial complex $X$, the $k$-skeleton of $X$-denoted $\text{Skel}_{k}(X)$-is the subcomplex $\{\alpha \in X: \dim (\overline{\alpha}) \leq k\}$. A homology $d$-sphere $S$ is said to be $k$-stacked if there is a homology $(d+1)$-ball $\Delta$ such that $S$ is the boundary of $\Delta$ and $\text{Skel}_{d-k} (\Delta) =\text{Skel}_{d-k}(S)$.

If $X$ is a $d$-dimensional simplicial complex then its face vector $(f_{-1}, f_0, \ldots, f_d)$ is defined by $f_i =f_i (X) =\# \{\alpha \in X: \dim (\overline{\alpha}) =i\}$, $-1 \leq i \leq d$. The $g$-vector $(g_0, g_1,\ldots,g_{d+1})$ of $X$ is defined by
\begin{eqnarray*}
g_j = g_j (X) =\sum\limits^{j-1}_{i=-1} (-1)^{j-i-1} \binom{d-i+1}{j-i-1} f_i (X), ~ 0 \leq j \leq d+1.
\end{eqnarray*}
These equations may be inverted to write the $f$-vector in terms of the $g$-vector.
\begin{eqnarray*}
f_i(X) =\sum\limits^{i+1}_{j=0} \binom{d-j+2}{i-j+1} g_j (X), -1 \leq i \leq d.
\end{eqnarray*}
Since each $f_i$ is a non-negative linear combination of the $g_j$'s, any lower bound on the $g_j$'s in terms of topological invariants of $X$ is called a generalized lower bound theorem for $X$. (The original lower bound theorem was the lower bound  $g_2 \geq 0$ for homology spheres.)

One reason for our singling out the tame spheres is the following observation. It says that the tame spheres satisfy the generalized lower bound conjecture (GLBC) of McMullen and Walkup [4], that too for an easy and intuitively obvious reason.

\begin{Theorem} (GLBT for tame spheres) Let $S$ be a tame $d$-sphere and $0 \leq \ell \leq d/2$. Then $g_{\ell+1}(S) \geq 0$. Also, $g_{\ell+1}(S)=0$ iff $S$ is $\ell$-stacked.
\end{Theorem}

Note that, by the Dehn-Sommerville equations, $g_{d/2+1} (S) =0$ for any homology sphere $S$ of even dimension $d$. Also, by [1, Theorem 2.9], all tame $d$-spheres are $\lceil \frac{d}{2} \rceil$-stacked.

Let $X$ be a simplicial complex of dimension $d$. Then the sigma-vector $(\sigma_0, \ldots,\sigma_d)$ of $X$ was defined in [2] by the formulae
\begin{eqnarray*}
\sigma_i =\sigma_i (X) =\sum\limits_{A \subseteq V(X)} \frac{\tilde \beta_i (X [A])}{\binom{f_0(X)}{\# (A)}}, ~ 0 \leq i \leq d.
\end{eqnarray*}
We also adopt the convention that $\sigma_i (X)=0$ when $i < 0$ or $i > \dim (X)$. In [2, Lemma 2.2] it was shown that the sigma-vector of any homology sphere $S$ of dimension $d \geq 2$ satisfies
\begin{eqnarray}
\left. \begin{array}{lll}
&& \sigma_d (S) =1,\\
&& \sigma_{d-1}(S) =1+\sigma_0 (S), \text{ and}\\
&& \sigma_{d-1-i} (S) =\sigma_i (S), 0 < i < d-1.
\end{array}
\right\}
\end{eqnarray}
We obtain the following precise estimate on (certain alternating sums of) the sigma numbers of a tame sphere. Compare [2, Theorem 3.5].

\begin{Theorem} Let $S$ be a tame $d$-sphere. Let $0 \leq \ell \leq d/2-1$ and $m=f_0(S)$. Then $\sum\limits_{i=0}^\ell (-1)^{\ell-i} \sigma_i (S) \leq \frac{m+1}{d+3} \sum\limits^{\ell+1}_{i=0} (-1)^{\ell +1-i} \frac{g_i (S)}{\binom{d+2}{i}}$. Equality holds here iff $S$ is $(\ell+1)$-stacked.
\end{Theorem}

We now advance:

\noindent{\bf Conjecture 1:} Theorem 1.3 is true of all homology $d$-spheres.

The two-dimensional case of Conjecture 1 was recently proved in [3, Theorem 1.1]. Thus we have:

\begin{Theorem} (Butler et al): Any $m$-vertex homology sphere $S$ of dimension 2 satisfies $\sigma_0(S) \leq \frac{m+1}{5} \sum\limits^1_{i=0} (-1)^{1-i} \frac{g_i (S)}{\binom{4}{i}}$. Equality holds here iff $S$ is 1-stacked.
\end{Theorem}
Of course, all homology 2-spheres are combinatorial (indeed polytopal) spheres.

The mu-vector of a simplicial complex was introduced in [2]. Here we modify the definition as follows.

\begin{Definition} (New mu) Let $X$ be a $d$-dimensional simplicial complex. The mu-vector $(\mu_0, \mu_1,\ldots,\mu_d)$ of $X$ is defined by the formulae
\begin{eqnarray*}
\mu_0 &=& \mu_0 (X) =\sum\limits_{x \in V(X)} \frac{1}{1+f_0 (\ell k(x,X))}, \text{ and }\\
\mu_i &=& \mu_i (X) =\sum\limits_{x\in V(X)} \frac{\delta_{i1} +\sigma_{i-1} (\ell k (x,X))}{1+f_0 (\ell k(x,X))}, ~ 1 \leq i \leq d.
\end{eqnarray*}
\end{Definition}

Here $\delta_{ij}$ is the usual Kronecker symbol: $\delta_{ij}=1$ if $i=j$ and $\delta_{ij}=0$ otherwise. Recall that a simplicial complex $X$ is said to be 2-neighbourly if any two vertices of $X$ form an edge. Clearly we have $f_0 (X) -f_0 (\ell k (x,X)) \geq 1$ for $x \in V(X)$. Equality holds here for all $x$ iff $X$ is 2-neighbourly. Therefore, Definition 1.5 implies that $\mu_0 (X) \geq 1$, with equality iff $X$ is 2-neighbourly. Indeed, this definition coincides with [2, Definition 2.1] precisely for 2-neighbourly simplicial complexes.

With the new definition, the mu-vector of a closed homology manifold still satisfies the duality of [2, Theorem 2.3]. Namely, we have,

\begin{Theorem} If $M$ is a closed homology $d$-manifold then $\mu_{d-i} (M) =\mu_i (M), ~ 0\leq i \leq d$.
\end{Theorem}

Moreover, [2, Theorem 2.6], which was proved with the hypothesis of 2-neighbourliness, is now valid for all simplicial complexes:

\newpage
\begin{Theorem} (Morse inequalities) Any simplicial complex $X$ of dimension $d$ satisfies:-
\begin{enumerate}
\item[(a)] For $0 \leq j \leq d, ~ \sum\limits_{i=0}^j (-1)^{j-i} \mu_i (X) \geq \sum\limits^j_{i=0} (-1)^{j-i} \beta_i (X)$, with equality for $j=d$.
\item[(b)] Equality holds in (a) for some index $j$ iff for every induced subcomplex $Y$ of $X$, the morphism $H_j (Y) \rightarrow H_j(X)$, induced by the inclusion map $Y \rightarrow X$, is injective.
\item[(c)] For $0 \leq j \leq d, ~ \mu_j (X) \geq \beta_j (X)$.
\item[(d)] Equality holds in (c) for some index $j$ iff for every induced subcomplex $Y$ of $X$, the two morphisms $H_{j-1}(Y) \rightarrow H_{j-1}(X)$ and $H_j(Y) \rightarrow H_j (X)$, induced by the inclusion map $Y \rightarrow X$, are both injective.
\end{enumerate}
Recall that a simplicial complex $X$ is said to be $\mathbb{F}$-tight if $X$ is connected and, for every induced subcomplex $Y$ of $X$, the morphism $H_{\ast} (Y) \rightarrow H_\ast (X)$, induced by the inclusion map $Y \rightarrow X$, is injective. As an immediate consequence of Theorem 1.7, we have
\end{Theorem}

\begin{Corollary} A $d$-dimensional simplicial complex $X$ is $\mathbb{F}$-tight iff $X$ is connected and $\mu_j (X) =\beta_j (X)$ for all $ 0\leq j \leq d$.
\end{Corollary}

In particular, if $X$ is $\mathbb{F}$-tight then $\mu_0 (X) =\beta_0 (X)=1$, and hence (by the remark following Definition 1.5) $X$ is 2-neighbourly. Since for 2-neighbourly simplicial complexes, the new definition of the mu vector agrees with the old one, in this paper we have nothing more to say about tightness.

Recall that a homology $d$-manifold is locally tame if all its vertex links are tame. Similarly, it is locally $\ell$-stacked if all its vertex links are $\ell$-stacked homology spheres. The following result is a straightforward consequence of Theorem 1.3. (Compare [2, Theorem 3.6].)

\begin{Theorem} Let $M$ be a locally tame $d$-manifold. Then, for $ 1\leq \ell \leq \frac{d-1}{2}$, $\sum\limits^\ell_{i=0} (-1)^{\ell-i} \mu_i (M) \leq (-1)^\ell +\frac{g_{\ell+1} (M)}{\binom{d+2}{\ell+1}}$. Equality holds here iff $M$ is locally $\ell$-stacked.
\end{Theorem}

We now pose:\\
{\bf Conjecture 2:} Theorem 1.9 is true of all closed homology $d$-manifolds.

Using Theorem 1.4, we are able to prove this conjecture in dimension 3. Thus, we have

\begin{Theorem} Any closed homology 3-manifold $M$ satisfies $\mu_1 (M)-\mu_0 (M) \leq -1 +\frac{1}{10}g_2 (M)$. Equality holds here iff $M$ is locally 1-stacked.
\end{Theorem}

In [5], the definition of stackedness was extended to arbitrary closed homology manifolds. Thus,
\begin{Definition} Let $0 \leq \ell \leq d$, and let $M$ be a closed homology $d$-manifold. We say that $M$ is $\ell$-stacked if there is a homology $(d+1)$-manifold $\Delta$ such that $\text{Skel}_{d-\ell} (\Delta)=\text{Skel}_{d-\ell} (M)$ and $M$ is the boundary of $\Delta$.
\end{Definition}

Obviously, every $\ell$-stacked manifold is also locally $\ell$-stacked, though the converse is false in general. As a consequence of Theorems 1.7 and 1.9, we get: (Compare [2, Theorem 3.7].)

\begin{Theorem} (GLBT for locally tame combinatorial manifolds) Let $M$ be a connected and locally tame $d$-manifold, and let $1 \leq \ell \leq \frac{d-1}{2}$. Then $g_{\ell+1} (M) \geq \binom{d+2}{\ell+1} \sum\limits^\ell_{i=1} (-1)^{\ell-i} \beta_i (M)$. Equality holds here iff $M$ is $\ell$-stacked.
\end{Theorem}

The last conjecture in this paper is the following (compare [2, Conjecture 1.6] and [5, Question 5.10]).

\noindent{\bf Conjecture 3:} (GLBC for closed homology manifolds). Theorem 1.12 is true of all connected and closed homology manifolds.

Actually, the ``if" part of conjecture 3 is an immediate consequence of Theorem 3.1 and Proposition 5.2 of [5]. Thus we have,

\begin{Theorem} (Murai and Nevo). Let $1 \leq \ell \leq \frac{d-1}{2}$ and let $M$ be an $\ell$-stacked connected and closed homology $d$-manifold. Then
\begin{eqnarray*}
g_{\ell+1} (M) = \binom{d+2}{\ell+1} \sum\limits^\ell_{i=1} (-1)^{\ell-i} \beta_i (M).
\end{eqnarray*}
\end{Theorem}

In case $\ell=1$, Conjecture 3 says that any connected and closed homology manifold $M$ of dimension $d \geq 3$ satisfies $g_2 (M) \geq \binom{d+2}{2}\beta_1 (M)$, with equality iff $M$ is 1-stacked. In [6, Theorem 5.2] Novik and Swartz proved this inequality under the assumption that $M$ is $\mathbb{F}$-orientable, i.e., $\beta_d (M) >0$. They also proved that, when $d\geq 4$, equality holds for $\mathbb{F}$-orientable closed homology manifolds $M$ iff $M$ is locally 1-stacked. (Actually, Novik and Swartz assumes that the field $\mathbb{F}$ is infinite. But this is not a serious restriction since $\beta_1 (M; \mathbb{F})$ depends only on the characteristic of $\mathbb{F}$. Also note that a closed homology manifold of dimension $d \geq 4$ is 1-stacked iff it is locally 1-stacked.) These authors asked if in the case $d=3$ also, equality in this inequality holds only if $M$ is locally 1-stacked. Using Theorem 1.10, we can prove that Conjecture 3 is valid in dimension three. Since 1-stacked manifolds are trivially locally 1-stacked, this answers the question of Novik and Swartz in the affirmative. We have:

\begin{Theorem} Let $M$ be a connected and closed homology 3-manifold. Then $g_2 (M) \geq 10 \beta_1 (M)$. Equality holds here iff $M$ is 1-stacked.
\end{Theorem}

Notice that we have no orientability hypothesis here. Of course, all homology 3-manifolds are actually combinatorial manifolds.

For any connected and closed topological 3-manifold $M$, Walkup [8] defined the topological invariant $\gamma(M)$ as the infimum of $f_1(X)-4f_0 (X)$ over all triangulations $X$ of $M$. He proved that $\gamma(M)$ is always finite, and computed it explicitly for a few 3-manifolds. Since $g_2(X) =f_1(X)-4f_0 (X)+10$ for 3-dimensional simplicial complexes $X$, Theorem 1.14 may be rephrased as follows.

\begin{Corollary} For any connected and closed topological 3-manifold $M$, $\gamma (M) \geq 10 (\beta_1 (M; \mathbb{F})-1)$. Equality holds here iff $M$ has an 1-stacked triangulation.
\end{Corollary}

This result indicates that the natural generalization of Walkup's invariant for connected topological closed manifolds $M$ of dimension $d$ are perhaps the numbers
\begin{eqnarray*}
\gamma_\ell (M) := \text{Inf} \{g_{\ell+1} (X) : X \text{ is a triangulation of } M\}, ~ 0 \leq \ell \leq \frac{d-1}{2}.
\end{eqnarray*}
Conjecture 3 would imply that, for those manifolds which can be triangulated, these invariants are finite. Perhaps it is easier to establish the finiteness of these invariants?

In [5], Murai and Nevo amplified the commutative algebraic arguments of Novik and Swartz [6,7] to show that when $M$ is an $\mathbb{F}$-orientable connected and closed homology manifold all whose vertex links have the weak Lefschetz property (WLP), the inequalities of Conjecture 3 hold for $M$; also, if $\ell < \frac{d-1}{2}$, equality holds precisely for $\ell$-stacked manifolds $M$. Conjecturally, all homology spheres have the WLP. Moreover, it may not be too hard to prove (perhaps by an induction on the length $\lambda(S)$, see Definition 2.2 below) that all tame spheres have the WLP. However, it appears that commutative algebra is unable to tackle the non-orientable case, or the case of equality when $\ell =\frac{d-1}{2}$. Investigation of the sigma vectors of general homology spheres may provide an alternative route to the GLBC for homology manifolds. Indeed, the arguments in this paper show that we have the following implications:

Conjecture 1 $\Rightarrow$ Conjecture 2 $\Rightarrow $ Conjecture 3.

\section{Proofs}

It is well known that the $g$-vector of a simplicial complex changes under a single bistellar move according to the following rule. If a simplicial complex $Y$ of dimension $d$ is obtained from a simplicial complex $X$ by a bistellar move of index $t$, then, for $0\leq j \leq d$,
\begin{eqnarray}
g_{j+1} (Y) -g_{j+1} (X) = \left\{ \begin{array}{lll} +1  \text{ if } j=t \neq d/2,\\
-1  \text{ if } j=d-t \neq d/2,\\
0  \text{ otherwise }.
\end{array}\right.
\end{eqnarray}
 Indeed, the simplicity of this rule may be the real reason why the $g$-vector (rather than the $f$-vector) is the right object to look at.

We shall have two occasions to use the following result. It may be of independent interest in the investigation of $\ell$-stacked spheres.

\begin{Lemma} Let $R$ and $S$ be homology spheres of dimension $d\geq 2\ell+1$. Suppose $R$ is $\ell$-stacked and $S$ is obtained from $R$ by a single bistellar move, say of index $t$. Then $S$ is $\ell$-stacked if and only if $t\neq \ell$.
\end{Lemma}

{\bf Proof:} Let $\alpha \mapsto \beta$ be the bistellar move of index $t$ by which $S$ is obtained from $R$. The homology sphere $S$ of dimension $d \geq 2\ell+1$ admits the $(d-t)$-move $\beta \mapsto \alpha$. If $S$ is $\ell$-stacked, then [1, Corollary 2.13] implies that $S$ does not admit any $(d-\ell)$-move; hence $t \neq \ell$. This proves the ``only if" part. Now suppose $t \neq \ell$. Then by [1, Corollary 2.13] applied to $R$, we get $t \leq \ell-1$ or $t \geq d-\ell+1$. Let $A$ be a homology $(d+1)$-ball such that $\partial A =R$ and $\text{Skel}_{d-\ell}(A)=\text{Skel}_{d-\ell}(R)$.

{\bf Case 1:} $0 \leq t \leq \ell-1$. Since $\beta \not\in R$, $\dim (\overline{\beta})=t\leq\ell-1 \leq d-\ell$ and $\text{Skel}_{d-\ell}(A) =\text{Skel}_{d-\ell}(R)$, it follows that $\beta \not\in A$. Since $\overline{\alpha} \ast \partial \beta \subseteq R \subseteq A$, we then have $A \cap (\overline{\alpha \cup \beta})=\overline{\alpha} \ast \partial \beta$. Put $B=A \cup \overline{\alpha \cup \beta}$. Since $A$ is a $(d+1)$-pseudomanifold, $\alpha \cup \beta$ is the only facet of $B$ not in $A$, and since $A\cap (\overline{\alpha \cup \beta})=\overline{\alpha} \ast \partial \beta \subseteq \partial A$, it follows that $B$ is also a $(d+1)$-pseudomanifold and $\partial B$ is obtained from $\partial A=R$ by the bistellar move $\alpha \mapsto \beta$. Thus, $\partial B=S$. Since the homology $(d+1)$-balls $A$ and $\overline{\alpha \cup \beta}$ meet in the homology $d$-ball $\overline{\alpha} \ast \partial \beta$ and the latter is contained in the boundary of both, Mayer-Vietoris theorem implies that $B$ is also a homology $(d+1)$-ball. Let $\gamma \in \text{Skel}_{d-\ell}(B)$. Since $\dim (\overline{\alpha})=d-t \geq d-\ell+1$, $\alpha \not\subseteq \gamma$. If $\gamma \in A$ then (as $A$ is $\ell$-stacked) $\gamma \in R$. Since $\gamma \in R$ and $\alpha \not\subseteq \gamma$, it follows that $\gamma \in S$. On the other hand if $\gamma \not\in A$ (but $\gamma \in B$) then ${\beta} \subseteq \gamma  \subset \alpha \cup \beta$ and hence $\gamma \in \overline{\beta} \ast \partial \alpha \subseteq S$. Thus $\gamma \in  S=\partial B$ in either case. So $\text{Skel}_{d-\ell}(B)=\text{Skel}_{d-\ell}(S)$ and $S=\partial B$. So $S$ is $\ell$-stacked.

{\bf Case 2:} $d-\ell+1 \leq t\leq d$. Let $\gamma \in \text{Skel}_{d-\ell} (\overline{\alpha \cup \beta})$. Since $\overline{\alpha} \ast \partial \beta \subseteq R \subseteq A$, if $\gamma \not\in A$ then it follows that $\gamma \supseteq \beta$ and hence $\dim (\overline{\gamma}) \geq \dim (\overline{\beta})=t\geq d-\ell+1$. This is a contradiction. So $\text{Skel}_{d-\ell} (\overline{\alpha \cup \beta}) \subseteq A$. If $\alpha \cup \beta \not\in A$, then a minimal non-face of $A$ contained in $\alpha \cup \beta$ would be a missing face of the $\ell$-stacked homology ball $A$, of dimension $\geq d-\ell+1 \geq \ell+1$, contradicting [1, Lemma 2.11]. Thus, $\alpha \cup \beta \in A$. Let $B$ be the pure simplicial complex whose facets are all the facets of $A$ except $\alpha \cup \beta$. Since $A$ is a $(d+1)$-pseudomanifold, so is $B$. Since $\alpha \cup \beta$ is a facet of $A$ and $(\partial A) \cap (\overline{\alpha \cup \beta})=\overline{\alpha} \ast \partial \beta$, it follows that $B\cap \overline{\alpha \cup \beta}=\overline{\beta} \ast \partial \alpha$, and $\partial B$ is obtained from $\partial A=R$ by the bistellar move $\alpha \mapsto \beta$. Thus, $\partial B=S$. Since the induced subcomplex of $\partial A$ on $\alpha \cup \beta$ is the $d$-ball $\overline{\alpha}\ast \partial \beta$, Mayer Vietoris theorem implies that, like $A,B$ is also a homology $(d+1)$-ball. Let $\gamma \in \text{Skel}_{d-\ell}(B)$. Since $\text{Skel}_{d-\ell}(B) \subseteq \text{Skel}_{d-\ell}(A) = \text{Skel}_{d-\ell} (R) \subseteq R$, it follows that $\gamma \in R$. Since $\gamma \in B$ and $\alpha \not\in B$, it follows that $\gamma \not\supseteq \alpha$. Since $\gamma \not\supseteq \alpha$ and $\gamma \in R$, we then have $\gamma \in S$. Thus $\text{Skel}_{d-\ell}(B)=\text{Skel}_{d-\ell}(S)$. Hence $S$ is $\ell$-stacked. $\hfill{\Box}$

We shall also use

\begin{Definition} Let $S$ be tame sphere of dimension $d$. Then the length $\lambda (S)$ of $S$ is the minimum number of bistellar moves of index $< d/2$ required to obtain $S$ from $S^d_{d+2}$.
\end{Definition}

Actually, the following argument shows that, if $S$ is a tame sphere of dimension $d$, then every sequence of bistellar moves of index $< d/2$, used to obtain $S$ from $S^d_{d+2}$, contains the same number of moves. Namely, we have
\begin{eqnarray*}
\lambda (S) =\sum\limits_{0 \leq \ell < \frac{d}{2}} g_{\ell+1}(S).
\end{eqnarray*}

{\bf Proof of Theorem 1.2:} Let $S$ be a tame $d$-sphere, and fix $0 \leq \ell \leq d/2$. The proof is by induction on $\lambda (S)$. If $\lambda(S)=0$ then $S=S^d_{d+2}$ satisfies $g_{\ell+1}(S)=0$, and $S$ is indeed $\ell$-stacked (being the boundary of the standard $(d+1)$-ball). So let $\lambda (S)  >0$. Then $S$ is obtained from a shorter tame $d$-sphere $R$ by a single bistellar move, say of index $t < d/2$. Equation (2) shows that we have $g_{\ell+1}(R) \leq g_{\ell+1}(S)$ (since $\ell \leq \frac{d}{2} < d-t$) with equality iff $t\neq \ell$. Since, inductively, $g_{\ell+1}(R) \geq 0$, it follows that $g_{\ell+1} (S)\geq 0$. Also, $g_{\ell+1} (S)=0$ iff $g_{\ell+1} (R)=0$ and $t \neq \ell$. Thus,  $g_{\ell+1}(S)=0$ iff $t \neq \ell$, and $R$ is $\ell$-stacked (by induction hypothesis). If $\ell =\frac{d}{2}$, then every tame $d$-sphere is $\ell$-stacked (by [1, Theorem 2.9]) and has $g_{\ell+1}=0$ (by Dehn-Sommerville equations). So let $\ell < \frac{d}{2}$. If $g_{\ell+1}(S)=0$, then $S$ is obtained from the $\ell$-stacked sphere $R$ by a move of index $t\neq \ell$, so that $S$ is $\ell$-stacked by Lemma 2.1. Conversely, suppose $S$ is $\ell$-stacked. Since $d-t > \frac{d}{2} > \ell$ and $R$ is obtained from the $\ell$-stacked sphere $S$ by a $(d-t)$-move, Lemma 2.1 implies that $R$ is $\ell$-stacked. So, by induction hypothesis, $g_{\ell+1}(R)=0$. Also, by Lemma 2.1, $t \neq \ell$. Therefore $g_{\ell+1}(S)=0$.
$\hfill{\Box}$

For any set $V$ and any non-negative integer $j$, we shall use the notation $\binom{V}{j}:= \{A \subseteq V: \# (A)=j\}$. In the following proofs, we also make repeated use of an well known identity for binomial coefficients. For any three non-negative integers $p,q,r$, we have
\begin{eqnarray}
\sum\limits^p_{i=0} \frac{\binom{p}{i}}{\binom{p+q+r}{r+i}} =\frac{p+q+r+1}{q+r+1} \binom{q+r}{r}^{-1}.
\end{eqnarray}
Our next result records the manner in which the sigma vector of a homology sphere changes under a single bistellar move. Unfortunately, its statement is complicated.

\begin{Lemma} Let $R$ and $S$ be homology spheres of dimension $d$. Suppose $S$ is obtained from $R$ by a single bistellar move, say of index $t$. Put $m=f_0(S)$ and set
\begin{eqnarray*}
c = \left\{ \begin{array}{lll} 1 & \text{ if }& 0 < t <d,\\
\frac{m+1}{m} & \text{ if } & t =0, \text{ and}\\
\frac{m+1}{m+2} &\text{ if } & t=d. \end{array} \right.
\end{eqnarray*}
Then, for $0 \leq i \leq d$, we have:
\begin{enumerate}
\item[\rm 1.] If $i \not\in \{t-1,t,d-t-1,d-t\}$ then $\sigma_i (S) =c \sigma_i(R)$.
\item[\rm 2.] Suppose $t \neq d/2$. Then,
\begin{enumerate}
\item[\rm (2a)] $\sigma_{t-1} (S) < c\sigma_{t-1} (R)$ (except when $t=0$), \\
$\sigma_{d-t-1} (S) > c \sigma_{d-t-1}(R)$ (except when $t=d$),\\
$\sigma_t (S) > c \sigma_t (R)$ and $\sigma_{d-t} (S) < c \sigma_{d-t}(R)$.
\item[\rm (2b)] $\sigma_t (S) -\sigma_{t-1} (S) =c (\sigma_t (R) -\sigma_{t-1}(R))+\frac{m+1}{d+3} \binom{d+2}{t+1}^{-1}$,\\
$\sigma_{d-t}(S)-\sigma_{d-t-1}(S) = c\left( \sigma_{d-t} (R) - \sigma_{d-t-1}(R)\right) -\frac{m+1}{d+3} \binom{d+2}{t+1}^{-1}$.
\end{enumerate}
\end{enumerate}
(Recall our convention: $\sigma_i =0$ for $i <0$. Therefore Part (2b) of this lemma actually gives the value of $\sigma_0 (S)$ in terms of $\sigma_0(R)$ when $t=0$ or $d$.)
\end{Lemma}

{\bf Proof:} Notice that, if we have the result for $t=0$, then we get the result for $t=d$ by interchanging $R$ and $S$ and replacing $m$ by $m+1$. Therefore, in the rest of the proof, we assume $0 \leq t < d$.

Let $U=V(R)$ and $V=V(S)$.

If $t=0$ then $U \multimap V$ (see Notation 2.5), and $A \mapsto A \cap U$ is a bijection between $\binom{V}{j}$ and $\binom{U}{j} \sqcup \binom{U}{j-1}$. Therefore we have:
\begin{eqnarray*}
\sum\limits_{j=0}^m \frac{1}{\binom{m}{j}} \sum\limits_{A \in \binom{V}{j}} \tilde \beta_i (R[A]) &=& \delta_{id} +\sum\limits^{m-1}_{j=0} \frac{1}{\binom{m}{j}} \sum\limits_{A \in \binom{V}{j}} \tilde \beta_i (R[A]) \\
&=& \delta_{id} +\sum\limits^{m-1}_{j=0} \frac{1}{\binom{m}{j}} \sum\limits_{B \in \binom{U}{j} \sqcup \binom{U}{j-1}} \tilde \beta_i (R[B])\\
&=& \sum\limits_{j=0}^{m-1} \left( \frac{1}{\binom{m}{j}} +\frac{1}{\binom{m}{j+1}} \right) \sum\limits_{B \in \binom{U}{j}} \tilde \beta_i (R[B]) \\
&=& \frac{m+1}{m} \sum\limits^{m-1}_{j=0} \frac{1}{\binom{m-1}{j}} \sum\limits_{B\in \binom{U}{j}} \tilde \beta_i (R[B])\\
&=& c \cdot \sigma_i (R).
\end{eqnarray*}
(Here, for the penultimate equality, we have used the trivial identity $\binom{m}{j}^{-1} +\binom{m}{j+1}^{-1} =\frac{m+1}{m} \cdot \binom{m-1}{j}^{-1}, ~ 0 \leq j \leq m-1$.)

Thus, when $t=0$, we have
\begin{eqnarray}
\sum\limits^m_{j=0} \frac{1}{\binom{m}{j}} \sum\limits_{A \in \binom{V}{j}} (\tilde \beta_i (S[A]) - \tilde \beta_i (R[A]) =\sigma_i(S)-c \sigma_i (R).
\end{eqnarray}
When $0 < t<d$, we have $c=1$ and $U=V$. Therefore, the equation (4) is actually valid for $0 \leq t <d$.

Now, let $\alpha \mapsto \beta$ be the bistellar move of index $t$ by which $S$ is obtained from $R$. Thus, $\alpha$ and $\beta$ are non-empty disjoint subsets of $V$ with $\# (\alpha)=d-t+1,\# (\beta)=t+1$. Let ${\cal A}$ (respectively ${\cal B}$) be the set of all subsets $A$ of $V$ such that $A \supseteq \beta$ and $A \cap \alpha =\varphi$ (respectively, $A \supseteq \alpha$ and $A \cap \beta=\varphi$). We also put:
\begin{eqnarray*}
{\cal A}^+ &=& \{A \in {\cal A} : \tilde \beta_t (S[A]) -\tilde \beta_t (R[A])=+1\},\\
{\cal A}^- &=& \{A \in {\cal A} : \tilde \beta_{t-1} (S [A])-\tilde\beta_{t-1} (R[A])=-1\},\\
{\cal B}^+ &=& \{A \in {\cal B} : \tilde \beta_{d-t-1} (S[A])-\tilde \beta_{d-t-1} (R[A])=+1\}, \text{ and }\\
{\cal B}^- &=& \{ A \in {\cal B} : \tilde \beta_{d-t} (S[A])-\tilde \beta_{d-t} (R[A])=-1\}.
\end{eqnarray*}
For any $A \subseteq V$, [2, Lemma 3.4] says that we have $\tilde \beta_i (S[A])=\tilde \beta_i (R[A])$, except when (i) $i=t, A \in {\cal A}^+$, or (ii) $i=t-1, A \in {\cal A}^-$, or (iii) $i=d-t-1, ~ A\in {\cal B}^+$ or (iv) $i=d-t, A \in {\cal B}^-$. Thus, Part 1 is immediate from equation (4). Also, when $i \in \{t-1,t,d-t-1,d-t\}$ and $t \neq d/2$, equation (4) simplifies as follows.
\begin{eqnarray}
\sigma_t (S) -c \sigma_t (R) &=& + \sum\limits^m_{j=0} \frac{1}{\binom{m}{j}} \sum\limits_{\underset{\# (A)=j}{A\in {\cal A}^+}} 1 \\
\sigma_{t-1} (S) -c \sigma_{t-1}(R) &=& -\sum\limits_{j=0}^m \frac{1}{\binom{m}{j}} \sum\limits_{\underset{\#(A)=j}{A \in {\cal A}^-}}1\\
\sigma_{d-t-1} (S) -c \sigma_{d-t-1} (R) &=& +\sum\limits^m_{j=0} \frac{1}{\binom{m}{j}} \sum\limits_{\underset{\#(A)=j}{A \in {\cal B}^+}} 1\\
\sigma_{d-t} (S) -c \sigma_{d-t} (R) &=& -\sum\limits^m_{j=0} \frac{1}{\binom{m}{j}} \sum\limits_{\underset{\#(A)=j}{A\in {\cal B}^-}} 1.
\end{eqnarray}
When $0 <t<d$ (so that $U=V$), we have $\beta \in {\cal A}^-, \alpha \in {\cal B}^+$ and (by Alexander duality and the exact sequence for pairs) $A \mapsto V\backslash A$ is a bijection between ${\cal A}^+$ and ${\cal B}^+$, as well as between ${\cal A}^-$ and ${\cal B}^-$. Thus, in this case, the sets ${\cal A}^{\pm}, {\cal B}^\pm$ are all non-empty. When $t=0$, we have ${\cal A}^-=\varphi$, but we still have $\alpha \in {\cal B}^+, U \in {\cal B}^-$ and $\beta \in {\cal A}^+$. Thus the sets ${\cal A}^\pm, {\cal B}^\pm$ are non-empty except that ${\cal A}^-=\varphi$ when $t=0$. Therefore, the right hand sides of equations (5) and (7) are strictly positive and those of equations (6) and (8) are strictly negative (excepting the RHS of (6) when $t=0$). This completes the proof of Part (2a).

The result [2, Lemma 3.4] also says that we have ${\cal A}={\cal A}^+ \sqcup {\cal A}^-$ and ${\cal B}={\cal B}^+ \sqcup {\cal B}^-$. Therefore subtracting equation (6) from (5) and (7) from (8), we get
\begin{eqnarray*}
\sigma_t(S)-\sigma_{t-1}(S) -c (\sigma_t (R) - \sigma_{t-1} (R)) &=& \sum\limits^m_{j=0} \frac{1}{\binom{m}{j}} \sum\limits_{\underset{\#(A)=j}{A \in {\cal A}}} 1 \\
&=& \sum\limits^m_{j=0} \frac{\binom{m-d-2}{j-t-1}}{\binom{m}{j}} \\
&=& \sum\limits^{m-d-2}_{i=0} \frac{\binom{m-d-2}{i}}{\binom{m}{i+t+1}} \\
&=& \frac{m+1}{d+3} \binom{d+2}{t+1}^{-1},
\end{eqnarray*}
(where in the last step we have used the identity (3)) and similarly $\sigma_{d-t} (S) -\sigma_{d-t-1}(S)-c(\sigma_{d-t}(R)-\sigma_{d-t-1}(R))=-\frac{m+1}{d+3} \binom{d+2}{t+1}^{-1}$. This completes the proof of Part (2b). $\hfill{\Box}$

{\bf Remark:} When $t=d/2$, the above argument gives
\begin{eqnarray*}
\sigma_{d/2} (S) -\sigma_{d/2} (R) = \sigma_{d/2-1} (S) -\sigma_{d/2-1}(R) =\sum\limits^m_{j=0} \frac{1}{\binom{m}{j}} \left( \sum\limits_{\underset{\# (A)=j}{A\in {\cal A}^+}} 1 -\sum\limits_{\underset{\#(A)=j}{A\in {\cal B}^-}} 1 \right).
\end{eqnarray*}
But we have no control over this expression.

As an immediate consequence of this lemma, we get:

\begin{Corollary} Let $R$ and $S$ be homology $d$-spheres. Suppose $S$ is obtained from $R$ by a bistellar move of index $t< \frac{d}{2}$. Let $m=f_0(S)$ and $0 \leq \ell \leq \frac{d}{2}-1$. Put $a_\ell (S)=\sum\limits_{i=0}^\ell (-1)^{\ell-i} \sigma_i (S)$, and define $a_\ell (R)$ similarly. Then we have:
\begin{itemize}
\item[(a)] If $t=0$, then $a_{\ell}(S)=\frac{m+1}{m} a_\ell (R)+(-1)^\ell \frac{m+1}{(d+2)(d+3)}$.
\item[(b)] If $1 \leq t \leq \ell$ then $a_\ell (S) =a_\ell (R)+(-1)^{\ell-t} \frac{m+1}{d+3}\binom{d+2}{t+1}^{-1}$.
\item[(c)] If $t \geq \ell+2$ then $a_\ell (S)=a_\ell (R)$.
\item[(d)] If $t=\ell+1$ then $a_\ell(S)< a_\ell(R)$.
\end{itemize}
\end{Corollary}

{\bf Proof of Theorem 1.3:} Let $S$ be a tame $d$-sphere, and fix $0 \leq \ell \leq d/2 -1$. The proof is by induction on $\lambda(S)$ (see Definition 2.2). If $\lambda(S)=0$, then $S=S^d_{d+2}$. In this case, $\sigma_i (S)=-\delta_{i0} ~ (0 \leq i < d)$, and $g_i(S)=\delta_{i0} ~ (0 \leq i \leq d+1)$. So equality holds. Also, $S$ is indeed $(\ell+1)$-stacked in this case.

So let $\lambda(S)>0$. Then, the tame $d$-sphere $S$ is obtained from a shorter tame $d$-sphere $R$ by a bistellar move of some index $t < \frac{d}{2}$. By induction hypothesis, the result is true of $R$.

If $t=0$, Corollary 2.4, induction hypothesis and equation (2) yield (in the notation of Corollary 2.4, which we continue to use):
\begin{eqnarray*}
a_\ell (S) &=& (-1)^\ell \frac{m+1}{(d+2)(d+3)}+\frac{m+1}{m} a_\ell(R) \\
&\leq& (-1)^\ell \frac{m+1}{(d+2)(d+3)} +\frac{m+1}{d+3} \sum\limits^{\ell+1}_{i=0} (-1)^{\ell+1-i} \cdot \frac{g_i (R)}{\binom{d+2}{i}} \\
&=& \frac{m+1}{d+3} \sum\limits^{\ell+1}_{i=0} (-1)^{\ell+1-i} \cdot \frac{g_i (R)+\delta_{i,1}}{\binom{d+2}{i}} \\
&=& \frac{m+1}{d+3}\sum\limits^{\ell+1}_{i=0} (-1)^{\ell+1-i} \frac{g_i(S)}{\binom{d+2}{i}}.
\end{eqnarray*}
Similarly, if $1 \leq t \leq \ell$, we have
\begin{eqnarray*}
a_\ell(S) &=& (-1)^{\ell-t} \frac{m+1}{d+3} \binom{d+2}{t+1}^{-1} +a_\ell (R)\\
&\leq& (-1)^{\ell-t} \frac{m+1}{d+3} \binom{d+2}{t+1}^{-1} + \frac{m+1}{d+3} \sum\limits^{\ell+1}_{i=0} (-1)^{\ell+1-i} \cdot \frac{g_i (R)}{\binom{d+2}{i}}\\
&=& \frac{m+1}{d+3} \sum\limits^{\ell+1}_{i=0} (-1)^{\ell+1-i} \cdot  \frac{g_i (R)+\delta_{i, t+1}}{\binom{d+2}{i}}\\
&=& \frac{m+1}{d+3} \sum\limits^{\ell=1}_{i=0} (-1)^{\ell+1-i} \cdot \frac{g_i (S)}{\binom{d+2}{i}}.
\end{eqnarray*}
If $\ell+2 \leq t < d/2$, then $a_\ell(S)=a_\ell(R)$ and $g_i(S)=g_i(R)$ for $0 \leq i \leq \ell+1$, so that the inequality for $S$ is equivalent to that for $R$, which we have by induction hypothesis.

If $t=\ell+1$, then again $g_i(S)=g_i(R)$ for $0 \leq i \leq \ell+1$, but now $a_\ell(S) < a_\ell (R)$, so that in this case the inequality for $R$ implies strict inequality for $S$. This completes the proof of the inequality.

The proof of the inequality also shows that equality holds for $S$ iff it holds for $R$ and $t\neq \ell+1$. By induction hypothesis, equality holds for $R$ iff $R$ is $(\ell+1)$-stacked. In case $\ell =\frac{d}{2}-1$, we have $t < \frac{d}{2}=\ell+1$, and both $R$ and $S$ are $(\ell+1)$-stacked (since by [1, Theorem 2.9], all tame spheres of dimension $2\ell+2$ are $(\ell+1)$-stacked). Hence equality holds for $S$ and $S$ is $(\ell+1)$-stacked in case $\ell=\frac{d}{2}-1$. So assume $\ell < \frac{d}{2}-1$. Since the reverse move from $S$ to $R$ has index $d-t > \frac{d}{2} \geq \ell+1$, Lemma 2.1 shows that $S$ is $(\ell+1)$-stacked iff $R$ is $(\ell+1)$-stacked and $t \neq \ell+1$. Therefore, $S$ satisfies equality iff $S$ is $(\ell+1)$-stacked. $\hfill{\Box}$

{\bf Proof of Theorem 1.6:} Let $0 < i < d-i <d$ (and hence $d\geq 3$). Applying equation (1) to the vertex links of $M$, we get
\begin{eqnarray*}
\mu_{d-i}(M) &=& \sum\limits_{x\in V(M)} \frac{\delta_{i,d-1}+\sigma_{d-i-1}(\ell k (x,M))}{1+f_0 (\ell k(x,M))}\\
&=& \sum\limits_{x\in V(M)} \frac{\delta_{i1}+\sigma_{i-1} (\ell k (x,M))}{1+f_0 (\ell k(x,M))}\\
&=& \mu_i (M).
\end{eqnarray*}
Also, for $d \geq 1$, $\sigma_{d-1}$ of any homology $(d-1)$-sphere is $1-\delta_{1d}$. Therefore,
\begin{eqnarray*}
\mu_d(M) &=& \sum\limits_{x \in V(M)} \frac{\delta_{1d}+\sigma_{d-1}(\ell k(x,M))}{1+f_0 (\ell k(x,M))}\\
&=& \sum\limits_{x\in V(M)} \frac{1}{1+f_0 (\ell k(x,M))}\\
&=& \mu_0 (M).
\end{eqnarray*}
$\hfill{\Box}$

To state the next lemma, we need:

\begin{Notation}
We use the notation $\multimap$ for the covering relation for set inclusion. Thus, for sets $A,B, A \multimap B$ iff $A \subseteq B$ and $\# (B \backslash A)=1$.
\end{Notation}
The following lemma is just Lemma 2.4 of [2], except that, with our new definition (Definition 1.5) for the mu-vector, we do not need the hypothesis of 2-neighbourliness.

\begin{Lemma} The mu-vector of any $m$-vertex $d$-dimensional simplicial complex $X$ is given by
\begin{eqnarray*}
\mu_i(X) = \frac{1}{m} \sum\limits^m_{j=1} \frac{1}{\binom{m-1}{j-1}} \sum\limits_{\underset{\underset{\# (B)=j}{A \multimap B}}{A,B \subseteq V(X)}} \beta_i (X[B],X[A]), ~ 0 \leq i \leq d.
\end{eqnarray*}
\end{Lemma}

{\bf Proof:} Let $V=V(X)$. Also, for $x \in V$, let $L_x =\ell k(x,X),V_x =V(L_x), m_x=\# (V_x)=f_0(L_x)$. Notice that for $A \subseteq V \backslash \{x\}$, we have
\begin{eqnarray*}
\beta_i (X[A \sqcup \{x\}], X[A])=\beta_i (X[A_x \sqcup \{x\}], X[A_x])
\end{eqnarray*}
where $A_x = A \cap V_x$. Also, each set in $\binom{V_x}{k}$ has exactly $\binom{m-m_x-1}{j-k-1}$ pre-images in $\binom{V \backslash \{x\}}{j-1}$ under the map $A \mapsto A_x$. Therefore the right hand side (RHS) in the statement of this lemma equals
\begin{eqnarray*}
&& \frac{1}{m} \sum\limits^m_{j=1} \frac{1}{\binom{m-1}{j-1}} \sum\limits_{x\in V} \sum\limits_{A \in \binom{V \backslash x}{j-1}} \beta_i (X[A \sqcup x],X [A])\\
&=& \sum\limits_{x\in V} \sum\limits^{m_x}_{k=0} \left( \frac{1}{m} \sum\limits^m_{j=1} \frac{\binom{m-m_x-1}{j-k-1}}{\binom{m-1}{j-1}} \right) \sum\limits_{A \in \binom{V_x}{k}} \beta_i (X[A \sqcup x],X[A]).
\end{eqnarray*}
But, for $0 \leq k \leq m_x \leq m-1$, the identity (3) yields
\begin{eqnarray*}
\frac{1}{m} \sum\limits^m_{j=1} \frac{\binom{m-m_x-1}{j-k-1}}{\binom{m-1}{j-1}} =\frac{1}{m} \sum\limits_{i=0}^{m-m_x-1} \frac{\binom{m-m_x-1}{i}}{\binom{m-1}{i+k}} =\frac{1}{(m_x+1)\binom{m_x}{k}}.
\end{eqnarray*}
Therefore,
\begin{eqnarray*}
\text{RHS } = \sum\limits_{x \in V} \frac{1}{m_x+1} \sum\limits^{m_x}_{k=0} \frac{1}{\binom{m_x}{k}} \sum\limits_{A \in \binom{V_x}{k}} \beta_i (X [A\sqcup x], X[A]).
\end{eqnarray*}
But, when $i=0, ~ \beta_i (X[A \sqcup x], ~ X[A])=\delta_{k, 0}$ where $k=\# (A)$. When $1 \leq i \leq d$, the excision theorem and the exact sequence for pairs yield
\begin{eqnarray*}
\beta_i (X[A\sqcup x], X[A]) =\left\{ \begin{array}{lll} \tilde \beta_{i-1} (L_x [A]) \text{ if } A \neq \varphi\\
\tilde \beta_{i-1} (L_x [A]) +\delta_{i1} \text{ if } A=\varphi. \end{array}\right.
\end{eqnarray*}
Therefore, we have the result in both cases. $\hfill{\Box}$

{\bf Proof of Theorem 1.7:} Note that the statement of this theorem is identical with that of [2, Theorem 2.6] except that (we use the new definition of the mu-vector and) the hypothesis of 2-neighbourliness is missing. Its proof is identical with the proof in [2] of the latter theorem, except that we appeal to Lemma 2.6 above in place of [2, Lemma 2.4]. $\hfill{\Box}$

It is well known (and easy to prove) that, for any simplicial complex $X$ of dimension $d$, we have
\begin{eqnarray}
\sum\limits_{x \in V(X)} g_i (\ell k(x,X))=(d+2-i) g_i(X) +(i+1)g_{i+1} (X)
\end{eqnarray}
for $0 \leq i \leq d$.

{\bf Proof of Theorem 1.9:} By hypothesis, for each $x\in V(M)$, $\ell k(x,M)$ is a tame $(d-1)$-sphere, and $0 \leq \ell -1 \leq \frac{d-1}{2}-1$. Therefore, applying Theorem 1.3 to $\ell k (x,M)$ (with $d-1$ and $\ell-1$ in place of $d$ and $\ell$), we get
\begin{eqnarray*}
\sum\limits^\ell_{i=1} (-1)^{\ell-i} \frac{\sigma_{i-1} (\ell k (x,M))}{1+f_0 (\ell k(x,M))} \leq \frac{1}{d+2} \sum\limits^\ell_{i=0} (-1)^{\ell-i} \frac{g_i (\ell k (x,M))}{\binom{d+1}{i}}.
\end{eqnarray*}
Adding this inequality over all $x\in V(M)$ and using (9), we find
\begin{eqnarray*}
\sum\limits^\ell_{i=0} (-1)^{\ell-i} \mu_i (M) &\leq& \frac{1}{d+2} \sum\limits^\ell_{i=0} (-1)^{\ell-i} \frac{(d+2-i) g_i(M)+(i+1)g_{i+1}(M)}{\binom{d+1}{i}}\\
&=& \sum\limits^\ell_{i=0} (-1)^{\ell-i} \left( \frac{g_i(M)}{\binom{d+2}{i}}+ \frac{g_{i+1}(M)}{\binom{d+2}{i+1}}\right)\\
&=& (-1)^\ell +\frac{g_{\ell+1}(M)}{\binom{d+2}{\ell+1}}.
\end{eqnarray*}
Thus we have the inequality. (The change in the range of summation in the leftmost sum is not a mistake!). Equality holds here iff each vertex link of $M$ satisfies equality in Theorem 1.3, and - by Theorem 1.3 itself - this holds iff each vertex link is $\ell$-stacked. $\hfill{\Box}$

{\bf Proof of Theorem 1.10:} The closed 3-manifold $M$ need not be locally tame. So we apply Theorem 1.4 (in place of Theorem 1.3) to the vertex links of $M$, and argue exactly as in the proof of Theorem 1.9 (with $d=3$, $\ell=1$). $\hfill{\Box}$

{\bf Proof of Theorem 1.12:} By assumption, $M$ is connected. Thus $\beta_0(M)=1$. By Theorems 1.7 and 1.9, we get $\sum\limits_0^\ell (-1)^{\ell-i} \beta_i (M) \leq (-1)^\ell +\frac{g_{\ell+1}(M)}{\binom{d+2}{\ell+1}}$. Cancelling $(-1)^\ell \beta_0 (M) =(-1)^\ell$ from both sides, we have the required inequality. By Theorems 1.7 and 1.9, equality holds here iff $M$ is locally $\ell$-stacked and, for every induced subcomplex $Y$ of $M$, the morphism $H_\ell (Y) \rightarrow H_\ell(M)$ is injective.

If $M$ is $\ell$-stacked, then equality holds by Theorem 1.13. Conversely, suppose $M$ satisfies equality. Since for $\ell < \frac{d-1}{2}$ every locally $\ell$-stacked closed homology $d$-manifold is actually $\ell$-stacked by [1, Theorem 2.20]=[5, Theorem 4.6], we are done in this case. So, let $d=2\ell+1$. Let's put $\Delta = \{\alpha \subseteq V(M) : \text{Skel}_{\ell+1} (\overline{\alpha})\subseteq M\}$. That is, $\Delta$ is the largest simplicial complex (with respect to set inclusion) whose $(\ell+1)$-skeleton agrees with that of $M$. To prove that $M$ is $\ell$-stacked, it suffices to show that $\Delta$ is a homology $(2\ell+2)$-manifold whose boundary is $M$.

For $x\in V(M)$, let $L_x =\ell k (x,M)$ and put $\overline{L}_x =\{\alpha \subseteq V(L_x): \text{Skel}_\ell (\overline{\alpha}) \subseteq L_x\}$. Since each $L_x$ is an $\ell$-stacked homology sphere of dimension $2\ell$, [1, Theorem 2.12] shows that $\overline{L}_x$ is a homology $(2\ell+1)$-ball with $L_x$ as its boundary.

{\bf Claim:} $\ell k (x,\Delta)=\overline{L}_x, ~ x \in V(\Delta) =V(M)$.

Indeed, if $\alpha \in \ell k (x,\Delta)$ then $\text{Skel}_{\ell+1} (\alpha \sqcup x) \subseteq M$ and hence $\text{Skel}_\ell (\alpha) \subseteq L_x$. Thus, $\ell k (x,\Delta) \subseteq \overline{L}_x$. To prove the claim, we need to show that, conversely, $\alpha \in \overline{L}_x \Rightarrow \alpha \sqcup x \in \Delta$, that is, $\text{Skel}_\ell (\overline{\alpha}) \subseteq L_x \Rightarrow \text{Skel}_{\ell+1} (\alpha \sqcup x) \subseteq M$. Clearly, it suffices to prove this in case $\dim (\overline{\alpha})=\ell+1$ (then the general case follows). So suppose, if possible, that $\dim (\overline{\alpha})=\ell+1$, $\text{Skel}_\ell (\overline{\alpha}) \subseteq L_x$ but $\text{Skel}_{\ell+1} (\alpha \sqcup x) \not\subseteq M$, i.e., $\alpha \not\in M$. Then $x \ast \partial \alpha$ is an induced subcomplex of $M$. Since the composed morphism $H_\ell (\partial \alpha) \rightarrow H_\ell (x \ast \partial \alpha) \rightarrow H_\ell (M)$ is injective, it follows that the morphism $H_\ell (\partial \alpha) \rightarrow H_\ell (x \ast \partial \alpha)$ is injective. But this is absurd since $H_\ell (\partial \alpha)=\mathbb{F}$ and $H_\ell (x\ast \partial \alpha)=\{0\}$. This contradiction proves the claim.

By the claim, all the vertex links of $\Delta$ are homology balls of dimension $2\ell+1$. So $\Delta$ is a homology manifold of dimension $2\ell+2$. Also, for each $x\in V(\Delta) =V(M)$, the claim gives $\ell k (x,\partial \Delta)=\partial \ell k (x,\Delta)=\ell k (x,M)$. Therefore $\partial \Delta =M$. Thus, $M$ is $\ell$-stacked in the case of equality. $\hfill{\Box}$

{\bf Proof of Theorem 1.14:} This proof is identical with (the case $d=3, \ell =1$ of) the proof of Theorem 1.12, except that we appeal to Theorem 1.10 (which does not require the hypothesis of local tameness) in place of Theorem 1.9. $\hfill{\Box}$

\end{document}